# EXPLORING STUDENTS' UNDERSTANDING OF MEDIAN

| | |
|---|---|
| Md Amiruzzaman | Karl W. Kosko |
| Kent State University | Kent State University |
| mamiruzz@kent.edu | kkosko1@kent.edu |

*This paper explores middle-grade students' conceptions of median. Describes, where and why they struggle and provides learning trajectory to improve their understanding.*

## Background and Overview

Measures of Central Tendency (MOCT), including mean, mode and median, are important mathematical concepts typically included in middle grades standards and curriculum (National Council of Teachers of Mathematics [NCTM], 2000). However, the vast majority of research in MOCT has generally focused on the teaching and learning of arithmetic mean, with little study of median or mode (Groth & Bergner, 2006; Mokros & Russell, 1995). The few studies focusing on students' conceptions of median have found that high school and college-level students tend to have more difficulty calculating median than arithmetic mean (Barr, 1980; Zawojewski & Shaughnessy, 2000). Barr (1980) studied high school students' knowledge of median and found that less than half of participating students calculated median correctly and tended to provide incorrect explanations for their answers. Specifically, many high school students did not grasp the concept of median as the middle value in a distribution of data. Similarly, Zawojewski & Shaughnessy (2000) observed that although median is a simpler arithmetic calculation as compared to mean, most undergraduate students in their sample had more difficulty calculating median than statistical mean. Although the literature on students' understanding of median suggests that students have more difficulty with this concept than arithmetic mean, such research does not focus on the nature of such conceptions, or how and why they develop in students, indicating a significant gap in the field's models of students' learning and conceptions of median. Therefore, the purpose of the present study is to explore and model two middle grade students' conceptions of median. The study focused on a small sample of middle grades students both to allow for more in depth analysis of how conceptions are conveyed through student actions and because the middle grades are typically the initial point in K-12 curricula where these concepts are formally studied (NCTM, 2000).

## Method

We used a classroom based teaching experiment to build our models of middle grade students' conceptions of median and mode. A *classroom based teaching experiment* involves the teacher-researcher (the first author) in the posing of tasks to construct models of students' conceptions, but also to provide tasks that could scaffold students' learning. As with individual teaching experiments, the approach requires the teacher-researcher to collaborate with an observer (Behr, Wachsmuth, Post, & Lesh, 1984). The approach allows the teacher-researcher to interview students individually or as a group (Behr et al., 1984), and extends the individual-based constructivist teaching experiment to the social context.

The study was conducted in an urban, middle-class school district with 6 sixth-grade students (5 males and 1 female students). Students participated in 14 teaching experiment episodes across 14 weeks in Spring 2015. All episodes were video and audio recorded, with the observer-researcher taking field notes. Each episode took place in a designated room, isolated from external distractions and noise, and all participants sat with the teacher-researcher at a round table. This environment provided the necessary proximity to establish teacher-student relationships and observe students' actions more easily. Episodes involved exchanges in which the teacher-researcher presented a task focusing on MOCT (mean, median, and/or mode). Participants would work on the task either individually or with peers, and the teacher-researcher would ask probing questions to solicit descriptions of their mathematical thinking. However, each student was asked to explain their answers separately. After each student explained their answers, their written work was collected for further data analysis.

Each episode lasted approximately 45-50 minutes. Tasks for initial episodes focusing on a topic (i.e., median) were prepared ahead of time. However, tasks for all other episode were prepared based on the analysis of the prior episode. Specifically, the teacher-researcher and observer-researcher met and discussed current observations and developed hypotheses following each episode. Based on the hypotheses, a new protocol, including tasks, was developed for the next episode to test the hypothesis(es). For sake of space and simplicity, the present study presents data from episodes 6 and 8 to illustrate preliminary findings for learning trajectories of the participating students Alice and Bob.

## Analysis and Findings

**Initial conceptions**

In episode 6, both Alice and Bob described the median as the middle value of a data set, and generally were able to identify the median with data sets with an odd number of elements. Figure 1 illustrates Alice's written strategy to find median from a data set with an odd number of elements. When asked to explain how she found the median, Alice replied, "Oh, it was easy…you know…you write the number small to big and cross out one number from left and one number from right until you find a number in the middle." This statement conveys a potential interpretation on Alice's part of the median as the middle value in a data set. However, when Alice was asked to find the median for data with an even number of elements, she claimed that there was no median. Figure 2 illustrates Alice's written strategy to find median from a data set with an even number of elements. Noticeably, 17 and 20 are left as is with no continuation additional work by Alice to find the median. As Alice explained, "there is no median because 2 numbers are in the middle…not one. Only one number can be in the middle, not two". Bob's work on both tasks was nearly identical to Alice's with corresponding explanations. Bob explained his position and said, "Median lives in the middle of all numbers…when there is 7 numbers then we have a middle number, but when we have 10 numbers then we do not have a middle number…no median".

*6,7,7,8,9,12,14*

**Figure 1.** Alice's response to data with an odd number of elements for median.

*4,4,4,17,17,20,24,27,29,34.*

**Figure 2.** Alice's response to data with an even number of elements for median.

Although a seemingly simple 'mis'conception of the definition of the median, we conjecture that both students' conception of median is more nuanced than a simple misinterpretation. Specifically, Alice's description of how to find the median focused on marking off the discrete elements in the data set (once they were sorted). At no point did Alice or Bob data set, but instead referred to the elements within the set. We conjecture that these references, both in the students' spoken descriptions and in their written work, indicate that Alice and Bob were considering the elements as discrete objects of the data set instead of seeing the elements as part of a unitized data set. Similar ways of operating mathematically have been described regarding children's counting schemes. Discussing Steffe's schemes for children's number sequences, Oliver (2000) described the difference between pre-numerical counting schemes in which the end number is not a unitized collection of 1s (i.e., a child may count 1,2,…5 but not view the 5 as five 1s), and initial number sequences where they are able to consider a number as a collection of 1s in activity. In the present study, Alice and Bob demonstrated a similar focus on operating on discrete elements in the data set to consider the median as an end result of those discrete elements, instead of as a property of the unitized data set. We refer to this initial scheme, demonstrated in episode 6 by both Alice and Bob as the middle-value scheme. *Middle-value* is the number that is in the middle of a data set with an odd number of elements.

**After the learning trajectory**

During the teaching experiment, we sought to promote perturbation in students' initial conceptions of median. Therefore, we asked students to begin plotting the data on a number-line. The number-line was chosen due to its potential for presenting the distribution of elements in a data set. Thus, we anticipated that the number-line might facilitate students' operating on data sets and not only data elements. So, students began using the number-lines to model their process for finding the median during episode 8. At one point in the episode, students were provided a data set with odd number of elements and asked to find the median. Alice and Bob initially wrote down each number on the number-line and tried to find which number lay in the center/middle (see Figure 3). Once they found the center/middle number they marked the number and claimed that number as median.

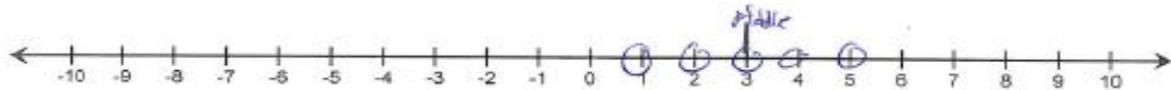

**Figure 3.** Alice's response to data with an odd number of elements for median.

After finding the median from a data set containing odd numbers of element, Alice and Bob were asked to find the median from a data set containing an even number of elements. Alice began to model the data on the number-line, as in the prior task, and used the number-line to sort the elements of the data on the number-line. Once all data were sorted, Alice identified the space between 2 and 3 (see Figure 4). Alice said, "this set has no middle value…but, we can find it…the middle number has to be in the middle of 2 and 3…we can add 2 and 3 and divide them in half…yes, then 2.5 is the median…(giggling)". Alice provided the response of 2.5 as a middle value of the given discrete elements. We interpret Alice' actions as an application of partitioning the data set when it was considered as a unitzed entity, or a *whole*. This is a distinctly different scheme than the middle-value scheme demonstrated in episode 6. We refer to this scheme as the center-value scheme because it involves operating on a unitized data set. *Center-value* is the calculated middle number from a data set with an even number of elements.

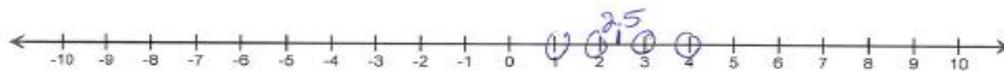

**Figure 4.** Alice's response to data with an even number of elements for median.

Similar to Alice, Bob said, "I see now…the middle number has to be between 2 and 3. I think, it is 2.5…it is the middle of 2 and 3…I am sure. Although it is tempting to consider the number-line as an intervention that allowed the student to clarify their definition of median, in episodes that followed, both Alice and Bob continued to use a number-line representation when asked to find median in future tasks. Specifically, when asked to find the median, students would sketch the portion of a number line representing the middle values in a sorted data set (often excluding all other values). This suggests a particular conception of median as the middle that has yet to be fully interiorized, but suggests a new scheme for operating on data to find median. Specifically, the student continued to determine median in activity and used the number line representation to scaffold such engagement. In episodes that followed, Alice and Bob were able to find the median. We conjecture this is due to the construction of a new scheme that allowed them to consider the data set as a whole (instead of only as discrete data elements). Therefore, they were more easily able to identify the middle position of the data set (using a segment of a number line) and determine the middle of the sorted data set as the median value.

## Conclusions

Findings of the present study indicate that middle grade students' conceptions of median may include middle-value and center-value schemes. Further, when presented a more visual representation of the data set, Alice and Bob both demonstrated construction of center-value schemes to consider the data set as a unitized whole, instead of as discrete data elements. Alice continued using a portion of a number line to find median, suggesting that the center-value scheme continued to be applied. These findings have several implications. A primary implication for research is that students' conception of the median is more complex than previously recognized. In regards to practice, and given the observed construction of the center-value scheme, these findings suggest that students should be afforded opportunities to manipulate data prior to being introduced to algorithms and purely symbolized definitions of these concepts.

## References

Barr, G. V. (1980). Some student ideas on the median and the mode. *Teaching Statistics*, *2*(2), 38-41.
Behr, M. J., Wachsmuth, I., Post, T. R., & Lesh, R. (1984). Order and equivalence of rational numbers: A clinical teaching experiment. *Journal for Research in Mathematics Education*, *15*(5), 323-341.
Groth, R. E., & Bergner, J. A. (2006). Preservice elementary teachers' conceptual and procedural knowledge of mean, median, and mode. *Mathematical Thinking and Learning*, *8*(1), 37-63.
Mokros, J. R., & Russell, S. J. (1995). Children's concepts of average and representativeness. *Journal for Research in Mathematics Education*, *26*(1), 20-39.
National Council of Teachers of Mathematics (NCTM). (2000). *Principles and standards for school mathematics: Discussion draft*. Reston, VA: Author
Zawojewski, J. S., & Shaughnessy, J. M. (2000). Mean and median: Are they really so easy? *Mathematics Teaching in the Middle School*, *5*(7), 436.